\documentclass[a4paper,UKenglish,cleveref, autoref, thm-restate]{lipics-v2021}
\pdfoutput=1
\authorrunning{T. Altenkirch and J. Neumann}

\usepackage{sty/lipics-packages}
\usepackage{sty/math-basic}
\usepackage{sty/include-macros}
\usepackage{sty/TT-meta-2}
\usepackage{sty/CwF-0}
\usepackage{sty/synSem}
\usepackage{sty/preload}

\theoremstyle{definition}

\newtheorem*{principle*}{Principle}

\newcommand{\ruleRef}[1]{(\nameref{tag:#1})}

\title{Synthetic 1-Categories in Directed Type Theory}
\titlerunning{Synthetic 1-Categories in Directed Type Theory}

\author{Thorsten Altenkirch}{University of Nottingham, United Kingdom}{thorsten.altenkirch@nottingham.ac.uk}{}{}
\author{Jacob Neumann}{University of Nottingham, United Kingdom}{jacob.neumann@nottingham.ac.uk}{https://orcid.org/0009-0003-5858-466X}{}

\authorrunning{T. Altenkirch and J. Neumann}

\Copyright{Thorsten Altenkirch and Jacob Neumann} 

\ccsdesc[500]{Theory of computation~Type theory}

\nolinenumbers 
\hideLIPIcs

\keywords{semantics, directed type theory, homotopy type theory, category theory,generalized algebraic theories}

\EventEditors{John Q. Open and Joan R. Access}
\EventNoEds{2}
\EventLongTitle{42nd Conference on Very Important Topics (CVIT 2016)}
\EventShortTitle{CVIT 2016}
\EventAcronym{CVIT}
\EventYear{2016}
\EventDate{December 24--27, 2016}
\EventLocation{Little Whinging, United Kingdom}
\EventLogo{}
\SeriesVolume{42}
\ArticleNo{23}
\begin{document}

\maketitle

\begin{abstract}
The field of \textit{directed type theory} seeks to design type theories capable of reasoning synthetically about (higher) categories, by generalizing the symmetric \textit{identity} types of Martin-L{\"o}f Type Theory to asymmetric \textit{hom-types}.
We articulate the directed type theory of the \textit{category model}, with appropriate modalities for keeping track of variances and a powerful directed-J rule capable of proving results about arbitrary terms of hom-types; we put this rule to use in making several constructions in synthetic 1-category theory.
Because this theory is expressed entirely in terms of \textit{generalized algebraic theories}, we know automatically that this directed type theory admits a syntax model and is the first step towards \textit{directed higher observational type theory}.
\end{abstract}

\section{Introduction}\label{introduction}

One exciting aspect of the emergent field of \textit{homotopy type theory (HoTT)}~\cite{hottbook} is the observation that \textit{types are $\infty$-groupoids}~\cite{van2011types}. Homotopy type theory can be understood as a \textit{synthetic} theory of $\infty$-groupoids: all the higher structure is generated by the simple rules for manipulating identity types in Martin-L{\"o}f Type Theory~\cite{MLTT1, MLTT2}, permitting efficient reasoning with these complex structures.

Not long after homotopy type theory was established, the search for \textit{directed homotopy type theory}---a synthetic theory of (higher) \textit{categories}---began. In a directed type theory, the identity types of ordinary Martin-L{\"o}f Type Theory (which are provably symmetric in the theory, i.e. a witness $p\colon\Id(t,t')$ can be turned into $p^{-1}\colon\Id(t',t)$) are replaced by asymmetric \textit{hom-types}. However, building a type theory to effectively work with these hom-types is beset by difficulties, in particular the need to carefully track the \textit{variances} of terms. A common feature of many approaches to directed type theory (e.g. \cite{licataHarper, nuyts, north2019}) is to track these variances by adopting some kind of modal typing discipline. However, no consensus ever emerged for exactly how to do this. More recent approaches to directed type theory (such as the work of Riehl and Shulman~\cite{riehl2017}) avoid these issues by adopting a more indirect approach inspired by simplicial spaces, at the cost of a more elaborate, multi-layered theory.

A possible new approach to directed type theory seems to be on the horizon, drawing from the recent development of \textit{higher observational type theory}~\cite{HOTT, altenkirch2024internal}. Higher observational type theory, or H.O.T.T., seeks to strike a balance between the properties of \enQuote{Book HoTT} (as originally articulated in \cite{hottbook}) and cubical type theory~\cite{cohen2018cubical}, particularly with regards to the central axiom of HoTT, Voevodsky's \textit{univalence axiom}. In Book HoTT, identity types are defined inductively, which results in the univalence axiom being impossible to compute with---complicating a key attribute of HoTT, its amenability to computer formalization. Cubical type theory rectifies this situation by instead defining identity types in terms of a special \textit{interval} type, and adding enough machinery to make univalence computable as a \textit{theorem} rather than an axiom (though at the expense of Book HoTT's intuitive simplicity). Higher observational type theory seeks to build homotopy type theory around \textit{observational} identities, in particular turning univalence into a \textit{definition}, preserving both the computational and intuitive character of homotopy type theory. H.O.T.T.'s definitional univalence suggests a clear directed analogue, making \textit{directed higher observational type theory} an appealing prospect. Though H.O.T.T. remains to be fully worked out, it's clear that \textit{second-order generalized algebraic theories (SOGATs)}~\cite{Uemura2023, uemura2021abstract, bocquet2022external, kaposi2024second} provide the appropriate setting for formulating this theory, as the language of SOGATs provide a \textit{higher-order abstract syntax} ideal for handling formal languages with variable binders.

Ordinary (first-order) generalized algebraic theories (GATs)~\cite{cartmell1986generalised} have played a prominent role in the semantics of type theory: Dybjer's \textit{categories with families (CwFs)}~\cite{dybjer} are a GAT articulating the basic mechanics of type theory, and provide a highly flexible and modular approach to interpreting type theories. Articulating the semantics of type theory as a GAT comes with numerous advantages: the universal algebraic features of GATs---such as homomorphisms, displayed models, products and coproducts of models, free and cofree models---are well-understood~\cite{kovacs2022type}. In particular, every GAT has an initial model, the \textit{syntax model}, which can be constructed as a quotient inductive-inductive type \cite{kaposi2019constructing}. Finally, GATs have the advantage of being (relatively) straightforward to formalize, as they make all the relevant operations and equations explicit.
Anticipating SOGAT and higher observational treatments of directed type theory, we begin by articulating directed type theory as a generalized algebraic theory.

\subsection{Related Work}
The present work draws most closely from Hofmann and Streicher's work on the groupoid model \cite{groupoidModel}; in particular, we develop a directed analogue of the groupoid model---the \textit{category model}---and adapt the groupoid models main constructs (dependent types and identity types) to the directed setting. We also closely follow the kind of metatheoretic arguments made there, and develop a directed analogue of their \textit{universe extensionality}, an early articulation of univalence. We also draw from the closely-related \textit{setoid model} of \cite{hofmann1995simple, altenkirch1999extensional}.

Among directed type theories, the present work draws some constructs from the theory of Licata and Harper \cite{licataHarper}, particularly their treatment of the \textit{opposite category} construction as a modality on contexts and context extension, as well as their treatment of $\Pi$-types in the directed setting. Our \textit{directed J-rule} for eliminating hom-types is similar to one of the eliminators given by North~\cite{north2019}, though with the critical difference mentioned below. Our theory, like North's, is \enQuote{1-dimensional} in the sense of Licata-Harper in that we maintain \textit{judgmental} equality as a symmetric notion, as opposed to \enQuote{2-dimensional} theories \cite{licataHarper, nuyts, Ahrens_2023} which introduce a theory of \textit{directed reductions}. All these theories adopt a modal typing discipline for handling variances, as do we, unlike the theories of \cite{riehl2017, kudasov2023formalizing, weinbergerModal} and \cite{weaverLicata}, which adopt approaches akin to simplicial and cubical type theories, respectively.

The commonality between the setoid, groupoid, preorder, and category model described in \autoref{catModel} is made much more general \cite{kovacs2022type, kaposi2019constructing}, where it is shown that any GAT gives rise to a CwF of algebras, algebra morphisms, displayed algebras, and sections. The \textit{modal} features of the present work are most likely instances of the modal type theory of \cite{gratzer2020multimodal}.

\subsection{Contribution and Organization}
We articulate a directed type theory satisfying the following constraints.
\begin{enumerate}
\item[(1)] It is presented as a generalized algebraic theory.
\item[(2)] It is 1-dimensional in the sense of Licata-Harper: there are no \enquote{directed reductions} introduced judgmentally.
\item[(3)] It is \textit{deeply-polarized}: there is a modal typing discipline to keep track of variances, which operates not just on types but on contexts, substitutions, and context extension.
\item[(4)] The directed J-rule (directed path induction) permits reasoning about arbitrary terms of hom-types.
\item[(5)] Hom-types can be iterated,\footnote{In contrast to e.g. \cite{licataHarper}, where homs-between-homs and homs-between-homs-between-homs is not possible to express.} expressing synthetic higher categorical structure (though in the present work we only consider 1-category theoretic structure).
\end{enumerate}
To our knowledge, there is no existing type theory satisfying all these criteria.

Starting with \autoref{catModel}, we adopt a semantics-driven approach by investigating a particular model, the \textit{category model} and abstracting its key features into a series of abstract notions of model (\autoref{DCwFs}). These notions are all GATs (indeed, CwFs with additional structure), and therefore each give rise to a syntax model. Our main notion is that of a \textit{Directed CwF (DCwF)}, a generalized algebraic theory of directed types with adequate polarity structure to properly track variances.

Achieving (4) while maintaining a modal typing discipline requires a novel approach. In the typing rules of existing directed type theories (including ours), the endpoint terms $t$ and $t'$ of a hom-type $\hom(t,t')$ are assigned opposite variances: $t$ negative and $t'$ positive. However, this poses a difficulty for typing the identity morphism $\refl_t\colon\hom(t,t)$ since $t$ must assume both variances. North~\cite{north2019} solves this by restricting $t$ to be a term of a \textit{core type} (interpreted semantically by groupoids), but the consequent J-rule only operates on hom-terms with a core endpoint, not arbitrary ones. Our solution instead uses groupoid \textit{contexts} rather than groupoid \textit{types}.

In \autoref{syntheticCT}, we show that this is a viable framework for conducting synthetic category theory. In this section, we adopt an informal style reminiscent of \cite{hottbook}, showing how this theory can be operated and how the groupoid context can be carefully maintained by a simple syntactic rule. We use our directed J-rule to give several basic constructions in synthetic (1-)category theory.

Finally, we consider the directed universe of sets in the category model, which serves as the category of sets. The existence of a directed universe allows us to make the metatheoretic argument that the syntax of DCwFs cannot prove the symmetry of hom-types (i.e. this is a genuinely \textit{directed} type theory) or the uniqueness of homs (analogous to Hofmann and Streicher's proof that the groupoid model refutes the uniqueness of identity proofs). We conclude by sketching several possible routes for further study.

\subsection{Metatheory and Notation}
Throughout, we work in an informal type-theoretic metatheory, using pseudo-\textsc{Agda} notation to specify GATs, make category-theoretic constructions, and define terms in the syntax of Directed CwFs. We use the notations
\[ (x\colon X)\to P(x) \qtq{and} (x\colon X)\times P(x) \]
for the dependent function and dependent sum types, respectively. When defining dependent functions, we'll enclose arguments in curly brackets to indicate that they're implicit. Any variables appearing free are also assumed to be implicitly universally quantified. We sometimes use underscores to indicate where the arguments to a function are written. When defining an instance $T$ of a construct given as a \code{record} type, we'll often omit the names of specific components, referring to all of them as just $T$ (matching the category-theoretic convention of referring to both the object- and morphism-parts of a functor $F$ by just $F$).

We use $=$ to mean \textit{definitional} or \textit{judgmental} equality in our metatheory, whereas $\equiv$ means \textit{propositional} equality (though there's no reason they couldn't coincide, i.e. in an extensional metatheory). We tacitly make use of appropriate extensionality principles (particularly function extensionality) for both notions of equality, and the uniqueness of identity proofs for $\equiv$. We write $\PROP$ for the type of h-propositions in our metatheory, i.e. those $P$ such that $p\equiv p'$ for all $p,p'\colon P$.

We assume basic familiarity with category theory. The set of objects of a category $\Gamma$ is denoted $\abs{\Gamma}$, the set of $\Gamma$-morphisms from $\gamma_0$ to $\gamma_1$ is denoted $\Gamma\:[\gamma_0,\gamma_1]$, and identities are written as $\iden$.  The \textit{discrete} groupoid/category on a set $X$ is the category whose objects are elements of $X$ and whose morphisms from $x_0$ to $x_1$ are inhabitants of the identity type $x_0\equiv x_1$. The \textit{opposite category} construction is understood to be definitionally involutive, i.e. $\abs{\Gamma\op}$ is defined to be $\abs{\Gamma}$ and $\Gamma\op\:[\gamma_0,\gamma_1]$ is defined to be $\Gamma\:[\gamma_1,\gamma_0]$, and thus
\[  (\Gamma\op)\op = \Gamma. \]

\section{The Category Interpretation of Type Theory}\label{catModel}
As mentioned, generalized algebraic theories (GATs) are a desirable formalism for expressing models of type theory, particularly when modelling numerous extensions to a \enquote{basic} type theory. When a theory is given as a GAT, all operations and equations are made clear and explicit, making it easier to compare and contrast similar theories. The theory of \keyword{Categories with Families (CwFs)}  (originally defined by Dybjer~\cite{dybjer}) present the fundamental operations of type theory---contexts, variables, terms, types, and substitutions---encoded as a GAT; upon this basic framework, an endless variety of different type theories can be studied.

\begin{FigCT}{CwF}{CwFDefn}

\caption{Main components of a CwF}
\end{FigCT}

The main components of a CwF are given in \autoref{tag:CwFDefn}: a category $\Con$ of \textit{contexts}, whose morphisms are called \textit{substitutions}; a presheaf $\Ty$ on $\Con$ and a dependent presheaf $\Tm$ over $\Ty$; and a \textit{context extension} operation guaranteeing that $\Tm$ is \textit{locally representable}  (in the sense of \cite{bocquet2022external}).
The last line says that there is an isomorphism (natural in $\Delta$) between the type of pairs $(\sigma,t)$ with $\sigma\colon\Sub\;\Delta\;\Gamma$ and $t\colon\Tm(\Delta,A[\sigma])$ and the type of substitutions from $\Delta$ to $\Gamma\cExtend A$. The left-to-right direction of this isomorphism is denoted $\langle\underlines,\underlines\rangle$ and the opposite direction as $\pp\circ\underlines, \vv[\underlines]$, so $$ \tau\equiv\langle \pp\circ \tau, \vv[\tau]\rangle \qtq{and} \sigma\equiv\pp\circ\langle\sigma,t\rangle \qtq{and} t\equiv\vv[\langle\sigma,t\rangle] $$
for any $\sigma,t$ as above and $\tau\colon\Sub\;\Delta\;(\Gamma\cExtend A)$.

Two paradigm examples of CwFs are the \textit{setoid model} of \cite{hofmann1995simple, altenkirch1999extensional} and the \textit{groupoid model} of \cite{groupoidModel}. In the former, the contexts are setoids (i.e. sets equipped with equivalence relations), the types are families of setoids (functorially) indexed over their context setoid, and terms are given by the appropriate notion of \textit{section} of their type (see the \textsc{Coq} formalization of \cite{altenkirch2021constructing} for a precise definition). The groupoid model is quite similar: contexts are groupoids, types are families of groupoids functorially indexed over their context groupoid, and terms are the appropriate notion of section. Indeed, we can view the groupoid model as generalizing the setoid model: a setoid can be viewed as a groupoid whose hom-sets are \textit{subsingletons}  (or \textit{propositions}, in the terminology of \cite{hottbook}), sets with at most one element. In other words, the groupoid model is what results when the assumption of \textit{proof-irrelevance} is dropped from the setoid model.

Both these models provide interpretations for numerous type formers, in particular the dependent types, identity types, and universes characteristic of Martin-L{\"o}f Type Theory~\cite{MLTT1, MLTT2}. The difference in these models is reflected in the type theories they interpret: while both models permit arbitrary iteration of the identity type former (expressing identities between identities, and identities between identities between identities, and so on), these iterated identity types become trivial more quickly in the setoid model. More precisely, the setoid model validates the \textit{uniqueness of identity proofs} principle, meaning that any two terms of an identity type, $p,q\colon\Tm(\Gamma,\Id(x,y))$ are themselves identical, $\mathsf{UIP}(p,q)\colon\Tm(\Gamma,\Id(p,q))$. The groupoid model famously violates this principle: in the type theory of the groupoid model, there are types (in particular, the universe of sets) which are not \textit{h-sets}, i.e. they possess terms which are proved identical by multiple, \textit{distinct} identity proofs.

This provides a roadmap for how we might develop a model of directed type theory. Since directed type theory can be described as \enQuote{dependent type theory, but with \textit{asymmetric} identity types}, this leads us to suspect that models of directed type theory will result if we simply drop the assumption of symmetry from the setoid and groupoid models. A setoid without symmetry is a preorder, and a groupoid without symmetry is a category. A close inspection of the definition of the groupoid model reveals that nothing in its interpretation of \textit{just the CwF structure} requires symmetry (i.e. that morphisms are invertible), and thus we can define the preorder model of type theory and the \textbf{category model of type theory}. The category model is given by \autoref{tag:CatModelDefn}, minus the definition of context extension (which will be discussed more below). This is just a generalization of the groupoid model, obtained by dropping symmetry: contexts are categories (and substitutions are functors), types are families of categories, and terms are sections.
We wont focus on the preorder model here, but leave it to future work to develop the directed analogue of setoid-model-specific considerations. Instead, we'll highlight those features of the category model which are relevant for modelling directed type theory, before abstracting those features into the notion of a \textit{directed CwF} in the next section.

\begin{FigCT}{categoryModel/main}{CatModelDefn}

\caption{The CwF structure of the category model, excluding context extension.}
\end{FigCT}

While the basic CwF structure of the groupoid model doesn't require symmetry (i.e. that all morphisms are invertible), its interpretations of further type formers certainly do. After all, our hope is that by passing from the groupoid model to the category model, the symmetric identity types of the former will become asymmetric \textit{hom-types} in the latter. Consider the semantics of the identity type former in the groupoid model. Here, and henceforth, we define a type (in this case $\Id(t,t')$) by giving its object- and morphism-parts, which are denoted \code{obj} and \code{map} in \autoref{tag:CatModelDefn}, but here are both written as just $\Id(t,t')$.
\includeCT{groupoidModel/Id-former}
Here, the fact that $A(\gamma_1)$ is a \textit{groupoid} is used in an essential way (we must take the inverse of $t'(\gamma_{01})$), and hence this definition doesn't work in the category model. But notice the following: the term $t$ is in the \enQuote{negative} position (the domain) and the term $t'$ is in the \enQuote{positive} position. Fittingly, we only use the \textit{inverse} of $t(\gamma_{01})$---never $t(\gamma_{01})$ itself---and only use $t'(\gamma_{01})$ but not its inverse. This observation will provide the key to adapting this definition for the category model.

What is needed is for $t$ to be a \textit{contravariant} term of type $A$, while keeping $t'$ as \textit{covariant}. This difference can be articulated in the category model, using a fundamental construct from category theory: \textit{opposite categories}. A type $A\colon\Ty\;\Gamma$ in the category model consists of a family of categories $A(\gamma)$ for each object $\gamma\colon\abs{\Gamma}$ and a functor $A(\gamma_{01})\colon\CAT\:[A(\gamma_0),A(\gamma_1)]$ for each morphism $\gamma_{01}\colon\Gamma[\gamma_0,\gamma_1]$. Given such a family of categories $A$, we can form a new family $A^-$, where $A^-(\gamma)$ is defined as the opposite category of $A(\gamma)$. This extends to the morphism part as well, because any functor $f\colon\CAT\:[C,D]$ can be viewed as a functor on their opposites, $f\colon\CAT\:[C\op,D\op]$. Alternatively, we could view $A$ as a functor $\Gamma\to\CAT$, and define $A^-$ to be the composition of $A$ with the endofunctor $(\underlines)\op\colon\CAT\to\CAT$. We can state generally that the category model validates the following rule:
\[ \infer{A^-\colon\Ty\;\Gamma.}{A\colon\Ty\;\Gamma}  \]
If $t\colon\Tm(\Gamma, A^-)$, this means that the object part of $t$ will still send objects $\gamma\colon\abs{\Gamma}$ to objects of $A(\gamma)$, since $A(\gamma)$ and $A^-(\gamma)$ have the same objects. But observe the type of its morphism part: $$ t \colon (\gamma_{01}\colon\Gamma\:[\gamma_0,\gamma_1]) \to (A\;\gamma_1)\:[t\;\gamma_1, A\;\gamma_{01}\;(t\;\gamma_0)]. $$
This is precisely what we need to articulate the definition of hom-types in the category model: see \autoref{tag:CatModelHom}. This definition is almost exactly the same as the semantics of $\Id$ in the groupoid model, but with $t$ changed to be a term of $A^-$, thus eliminating the need for the categories $A(\gamma_i)$ to be groupoids. Here's the hom-type formation, expressed as a rule:
\[ \infer{\hom(t,t')\colon\Ty\;\Gamma.}{t\colon\Tm(\Gamma,A^-) & t'\colon\Tm(\Gamma,A)} \]
The type annotation of $t$ as a \enQuote{negative} term and the implicit annotation of $t'$ as \enQuote{positive} serve as a kind of \textit{modal typing discipline} for keeping track of the \textit{variances} of terms.

\begin{FigCT}{categoryModel/Hom}{CatModelHom}

\caption{Semantics of the $\hom$-type former in the category model}
\end{FigCT}

For now, we just state the formation rule for hom-types; introducing and eliminating terms of hom-types will require more machinery. To see what kind of machinery, let's instead consider dependent function types. Like with the formation of hom-types, $\Pi$-types involve positive and negative \enQuote{variance}: a function is contravariant in its argument and covariant in its result. Therefore, as we might expect, the interpretation of $\Pi$-types in the groupoid model (\cite[Section 4.6]{groupoidModel}) makes essential use of the invertability of morphisms in a groupoid. Again, it only comes into play when defining the \textit{morphism} part of the interpretation: the object part (reproduced in \autoref{tag:CatModelPi}) defines for each $\gamma\colon\abs{\Gamma}$ an auxiliary type $B_\gamma$ in context $A(\gamma)$, and then specifies the category $\Pi(A,B)\;\gamma$ with terms $\theta\colon\Tm(A(\gamma),B_\gamma)$ as objects. This works fine in the category model. However, defining the morphism part of $\Pi(A,B)$ requires a kind of negative variance \textit{deeper} than the shallow contravariance of $A^-$: in the type $\hom(t,t')$ it was a \textit{term} that occurred negatively ($t$), in the type $\Pi(A,B)$ it's a \textit{type} that occurs negatively.

To make sense of this, we must consider the \textit{opposite category} operation, not just as acting on each $A(\gamma)$ in a family of categories over a context $\Gamma$, but as acting on the contexts themselves. In the category model, we have the following rules.
\[
    \infer{\Gamma^-\colon\Con}{\Gamma\colon\Con}
    \qquad \infer{\sigma^-\colon\Sub\;\Delta^-\;\Gamma^-}{\sigma\colon\Sub\;\Delta\;\Gamma}
\]
That is, we can negate contexts and substitutions as well as types: $\Gamma^-$ is interpreted as $\Gamma\op$, and this operation is (\textit{co}variantly) lifted onto functors as before. Now consider the difference in the morphism parts of terms with these different kinds of variance.
\includeCT{categoryModel/morphParts}
This is why we referred to this as \enQuote{shallow} and \enQuote{deep} negation: the difference between the first two is that we've flipped around each $A(\gamma)$, whereas in the third term, the dependence of $A$ on $\Gamma$ has itself been flipped around ($A$ is now contravariant, so $A\;\gamma_{01}$ takes objects of $A(\gamma_1)$ to objects of $A(\gamma_0)$). It is this latter kind of contravariance that describes $A$'s position in $\Pi(A,B)$.

We need another ingredient to state the $\Pi$-type formation rule: \textit{negative context extension}. In the type $\Pi(A,B)$, $A$ appears negatively, i.e. we want $A$ to depend negatively on $\Gamma$, i.e. $A\colon\Ty(\Gamma^-)$; but $B$ appears positively---we want $B$ to depend \textit{covariantly} on $\Gamma$, plus a variable of type $A$. In the theory of CwFs, a type \enQuote{depending on a variable} of another type is encoded by context extension. In the category model, we in fact have \textit{two} context extension operations (see \autoref{tag:CatModelExtend}), corresponding to the two ways a type can depend on a context. The positive context extension operator, $\cExtend[+]$, obeys the usual isomorphism discussed above. For the negative extension, the isomorphism becomes:

\begin{equation} \label{tag:NegLocalRep}
(\sigma\colon\Sub\;\Delta\;\Gamma)\times (\Tm(\Delta^-,A[\sigma^-]^-)) \cong \Sub\;\Delta\;(\Gamma\cExtend[-]A)
\end{equation}
for any $A\colon\Ty\;\Gamma^-$. We'll write $\pair{\underlines,_-\underlines}$ for the left-to-right direction of this isomorphism, and write $\pp_{-,A}\colon\Sub\;(\Gamma\cExtend[-]A)\;\Gamma$ and $\vv_{-,A}\colon\Tm((\Gamma\cExtend[-]A)^-,A[\pp_{-,A}^-]^-)$ for the data obtained from applying the right-to-left direction to the identity morphism on $\Gamma\cExtend[-]A$.
With this, we have everything needed to give the semantics of the $\Pi$-type former; this is done in \autoref{tag:CatModelPi}. As expected for $\Pi$-types, we have an isomorphism between $\Tm(\Gamma,\Pi^+(A,B))$ and $\Tm(\Gamma\cExtend[-]A,B)$, the application and lambda-abstraction rules. This is omitted here for reasons of space, but included in the appendix with accompanying calculations; see \autoref{tag:CatModelPiFull}.

\begin{FigCT}{categoryModel/extend}{CatModelExtend}

\caption{Semantics of context extension in the category model}
\end{FigCT}

\begin{FigCT}{categoryModel/Pi}{CatModelPi}

\caption{Semantics of the $\Pi$-type former in the category model}
\end{FigCT}

Let's now return to the key feature of directed type theory, hom-types. Above, we gave just the formation rule for hom-types, but said nothing of how to introduce or eliminate terms of this type. Stating the introduction rule, the term \code{refl} inhabiting $\hom(t,t)$ for each term $t$ proves rather subtle. The difficulty stems from the mixed-variance problem mentioned in the introduction: since our formation rule demands the domain term $t$ be of type $A^-$ and the codomain term $t'$ to be of type $A$, it's not immediately clear how to make $\hom(t,t)$ well-formed. There is, in general, no way to coerce terms of type $A$ into terms of type $A^-$ or vice versa, and we have no rule permitting us to use a term in both variances.

As mentioned in the introduction, the solution to this problem presented in \cite{north2019} is to use \textit{core types}. This solution consists of asserting a new type $A^0$ for each $A$, equipped with coercions $\Tm(\Gamma,A^0) \to \Tm(\Gamma,A)$ and $\Tm(\Gamma,A^0)\to \Tm(\Gamma,A^-)$. Then, for a term $t\colon\Tm(\Gamma,A^0)$, it makes sense to write $\hom(t,t)$, as $t$ can be coerced to both the positive and negative modality, in order to fit the $\hom$ formation rule. From there, a directed J-rule can be stated for eliminating hom terms. The issue with this solution is that it forces homs to have core endpoints: the directed J-rule can \textit{only} be used to prove claims about homs anchored at a core term, and it's not clear that proofs about arbitrary homs can be made. For the synthetic category-theoretic claims we study below, this will prove to be an unacceptable restriction.

A solution which avoids this shortcoming is revealed by considering hom-types in the empty context. In the empty context, a type $A$ is the same thing as a category, and a term of type $A$ is the same thing as an object of type $A$ (we silently coerce between $\mathbbm{1}\to X$ and $X$). So then there's no difference between terms of type $A$ and terms of type $A^-$, since a category and its opposite have the same objects. Therefore, in the empty context, there is no mixed-variance problem, and we can state the introduction rule for $\refl_t$ simply by coercing $t$ to be positive and negative as needed.

This doesn't extend to arbitrary contexts: as we saw above, terms $t\colon\Tm(\Gamma,A^-)$ and $t'\colon\Tm(\Gamma,A)$ have different morphism parts. But here's the key observation: if $\Gamma$ is a \textit{groupoid}, then we can still coerce between $A$ and $A^-$: given $t\colon\Tm(\Gamma,A^-)$, we can obtain $-t\colon\Tm(\Gamma,A)$, and vice-versa. The definition is given in \autoref{tag:CatModelNeutCoe}; there (and henceforth), we use $\Gamma\colon\NCon$ to indicate that $\Gamma$ is a groupoid, and therefore can invert $\Gamma$-morphisms as needed. So, rather than introduce a new type $A^0$ whose terms can be either positive or negative, we have instead have identified those contexts---neutral contexts---where terms of the familiar types $A$ and $A^-$ can be inter-converted. Given this, we can introduce \code{refl}:
\[ \infer{\refl_t\colon\Tm(\Gamma,\hom(t,-t)).}{\Gamma\colon\NCon & A\colon\Ty\;\Gamma & t\colon\Tm(\Gamma,A^-)} \]
We only need to assert $\refl_t$ for $t$ of type $A^-$, because the analogous rule for $t'$ of type $A$ can be derived: given $t'\colon\Tm(\Gamma,A)$, we observe that $t'=-(-t')$, so $\refl_{-t'}\colon\Tm(\Gamma,\hom(-t',t'))$.

\begin{FigCT}{categoryModel/neut-coe}{CatModelNeutCoe}

\caption{Semantics of neutral-context coercion in the category model}
\end{FigCT}

Let's conclude this section by giving an eliminator for our hom-type, known as the \textit{directed J-rule} or \textit{directed path induction}. Following \cite[Section 4.10]{groupoidModel}, we study directed path induction in the empty context first, which can then be extended to an arbitrary neutral context. Given $A\colon\Ty\Empty$ and $t\colon\Tm(\Empty,A^-)$ and some $M\colon\Ty(\Empty\cExtend[+]A\cExtend[+]\hom(t,\vv))$, our goal is to be able to prove $M[t',p]$ for arbitrary $t'$ and $p$, just by supplying a term $m$ of $M[-t,\refl_t]$. Translated into the category model semantics: $A$ is a category, $t$ and $t'$ are objects of $A$, $M$ is a functor from the coslice category $t/A$ into $\CAT$, $p$ is an $A$-morphism from $t$ to $t'$, and $m$ is an object of the category $M(t,\iden_t)$. The key observation is that $p$ is then \textit{also} a morphism in the coslice category from $(t,\iden)$ to $(t',p)$.
\[
\begin{tikzcd}
    & t \arrow[swap]{dl}{\iden_t} \arrow{dr}{p} \\
    t \arrow[swap]{rr}{p} && t'.
\end{tikzcd}
\]
Therefore, $$M(p) \colon \CAT\:[M(t,\iden_t), M(t',p)]$$ and so the object part of this functor turns objects of $M(t,\iden_t)$ into objects of $M(t',p)$, that is, it turns terms $m\colon\Tm(\Empty,M[-t,\refl])$ into terms $$(\mathsf{J}_{t,M}\;m)\:[t',p] \colon \Tm(\Empty,M[t',p]).$$ And, since $M(\iden)$ is the identity functor, we have the $\beta$ law, saying that $\mathsf{J}_{t,M}\;m\:[-t,\refl_t]\equiv m$. The general law replaces $\Empty$ with an arbitrary neutral context:

\begin{equation} \label{tag:dirJ}
\infer{
\mathsf{J}_{t,M}\;m\;\colon\Tm(\Gamma\cExtend[+]A\cExtend[+]\hom(t[\pp_A],\vv),M)
}
{
    \begin{matrix}
        \Gamma\colon\NCon\\
        t\colon\Tm(\Gamma,A^-)\\
        m\colon\Tm(\Gamma,M[-t,\refl_t])
    \end{matrix}
    &
    \begin{matrix}
        A\colon\Ty\;\Gamma\\
        M\colon\Ty(\Gamma\cExtend[+]A\cExtend[+]\hom(t[\pp_A],\vv))\\
        \
    \end{matrix}
}
\end{equation}
but the category model interpretation---see \autoref{tag:CatModelJ}---essentially follows this same idea. If $M$ doesn't need to depend on the term of type $\hom(t,\vv)$, then we can instead use the simpler rule

\begin{equation} \label{tag:simpledirJ}
\infer{
\mathsf{J}_{t,M}\;m\;\colon\Tm(\Gamma\cExtend[+]A\cExtend[+]\hom(t[\pp_A],\vv),M[\pp_{\hom(t[\pp_A],\vv)}])
}
{
    \begin{matrix}
        \Gamma\colon\NCon\\
        t\colon\Tm(\Gamma,A^-)
    \end{matrix}
    &
    \begin{matrix}
        A\colon\Ty\;\Gamma\\
        M\colon\Ty(\Gamma\cExtend[+]A)
    \end{matrix}
    &
    \begin{matrix}
    \\
    m\colon\Tm(\Gamma,M[-t])
    \end{matrix}
}
\end{equation}
Note that the dependence on $\hom(t[\pp_A],\vv)$ is preserved in the conclusion, even if $M$ ignores it. In \autoref{syntheticCT} we put this rule to use in synthetic category theory constructions and proofs.

\begin{FigCT*}{categoryModel/J}{CatModelJ}

\caption{Semantics of directed path induction in the category model}
\end{FigCT*}

\section{Directed Categories with Families}\label{DCwFs}
The aim of the present work is not just to establish the category model as a suitable interpretation of directed type theory, but to abstract the category model to a general, abstract notion of \enquote{model} of directed type theory. Specifically, we wish to present this model notion as a generalized algebraic theory, that is, as a CwF with further structure.
We do so in several stages, progressively capturing more of the structure described in the previous section. In addition to making the complex and multifaceted notion of \enquote{directed CwF} more digestible, this approach will also give us several intermediate notions, each of which is worthy of further study in its own right. First, we encapsulate the \enquote{negation} structure.

\mkDefn[Polarized CwF]{PCwF}
So a PCwF is just a CwF equipped with context-, substitution-, and type-negation involutions. The fact that the type-negation operation is a natural transformation just says that it is stable under substitution, i.e. $A[\sigma]^- \equiv A[\sigma^-]$. Now, notably absent from this definition is the negative context extension operation $\cExtend[-]$; by this definition, a PCwF only has the positive one. This is because the negative operation is, in fact, definable: in the category model, the following equation holds for any $\Gamma$ and any $A\colon\Ty\;\Gamma$:
\begin{equation} \label{tag:ExtendMinusEquality}
(\Gamma\cExtend[+]A)^- = \Gamma^-\cExtend[-]A^-.
\end{equation}
Here we use the fact that $(\Gamma^-)^-=\Gamma$,\footnote{In the category model, this holds as a definitional equality, though in \defnRef{PCwF} we only asserted it propositionally, since we're defining a GAT.} and hence $A\colon\Ty\;(\Gamma^-)^-$, making the right-hand side well-formed. Consequently, we can turn this equation around to \textit{define} negative context extension: for $A\colon\Ty(\Gamma^-)$, let $\Gamma\cExtend[-]A$ be $(\Gamma^-\cExtend[+]A^-)^-$. The isomorphism characterizing $\cExtend[-]$ (\autoref{tag:NegLocalRep}) can then be proved as a consequence of the one for $\cExtend[+]$.

Also absent from \defnRef{PCwF} is any mechanism connecting the context/substitution negation endofunctor to the type-negation operation. It's unclear if this ought to be rectified, or if their connection is just a peculiarity of the category model. Not every CwF fits the same mold of \enQuote{contexts are structures, types are families of structures}, so it's not possible to require in general that the type-negation operation is just post-composition with the context-negation functor. There \textit{are} suitably abstract ways of connecting the two---for instance, we can note that the category model is \textit{democratic} in the sense of \cite[Defn. 3]{castellan2021categories} : there is an isomorphism $K$ between contexts $\Gamma$ and closed types; this isomorphism is compatible with both negation operations, in that $K(\Gamma^-)=K(\Gamma)^-$. However, we don't need need such strong assumptions for the results of \autoref{syntheticCT}, so we omit them from the general definition of PCwFs.

Of course, the category model and the preorder model are both examples of PCwFs, where the negation is the \enquote{opposite} construction. But so are the groupoid and setoid models. Indeed, the groupoid model is a \textit{sub-PCwF} of the category model: a groupoid $\Gamma$ is a category, and so it makes perfect sense to take the opposite category of $\Gamma$, obtaining $\Gamma^-$, which is also a groupoid. What makes the groupoid model a peculiar instance of a PCwF is that $\Gamma\cong\Gamma^-$ for every $\Gamma$, and $A(\gamma)\cong A^-(\gamma)$ for every $A$. It is what we'll call a \textit{symmetric PCwF}. The setoid model is also a symmetric PCwF, but strictly so: there, $\Gamma \equiv \Gamma^-$. The situation exemplified by the groupoid/category and setoid/preorder models---a symmetric sub-PCwF of another PCwF---is what we capture in our next notion.
\mkDefn[Neutral-Polarized CwF]{NPCwF}

In the third bullet point, note that $A$ is \textit{not} assumed to be in $\NTy\;\Gamma^-$; if it were, then $\Gamma^-\cExtend A\colon\NCon$ and then we could use the isomorphism $e$ for $\Gamma^-\cExtend A$ and the coercion operators to construct this. With $A$ being an arbitrary type, this is a genuine addition to the theory. The requirements of this point are somewhat \textit{ad-hoc}: these were just the principles needed in \autoref{syntheticCT} to be able to operate effectively with neutral contexts (and all are provable in the category model). Perhaps a more mature version of this theory will place more requirements on NPCwFs, but this is all the neutral-polar structure needed here.

Let us also note that a common feature in directed type theories (e.g. \cite{nuyts, north2019}) is to include \textit{core types}, i.e. an operation of the form $(\underlines)^0\colon \Ty\;\Gamma\to\NTy\;\Gamma$. In the category model, this is interpreted as applying the \textit{core groupoid} construction to each category $A(\gamma)$, producing a family of groupoids indexed over $\Gamma$. We might as well have a deep version too, operating on contexts $(\underlines)^0\colon\Con\to\NCon$. We won't need these features for the present work (and therefore don't endeavor to axiomatize them), but, once again, it's quite possible that it would be fruitful to add them in the future.

Recall that CwFs are not a single notion of model for a single type theory, but rather that CwFs encode the basic structural operations of type theory, upon which innumerable different type theories can be specified by defining the desired term- and type-formers. We have arrived at the same point in our development of a semantics for directed type theory: the notion of NPCwF consists solely of structural components, but nothing that actually allows for the construction of interesting types and terms. So let's rectify this by giving the directed analogue of the standard core of undirected type theory: identity types, dependent types, and universes.

We start with the directed analogue of identity types, hom-types.
\mkDefn[Directed CwF]{DCwF}
The naming of $\hom$ versus $\Id$ is suggestive: the types in a DCwF are supposed to function like synthetic categories (with $\hom$ encoding their morphisms), and the neutral types are synthetic groupoids, whose homs are symmetric like an identity type. This point is best illustrated by the following claim.
    \begin{proposition}\label{state:Id-symm}
Every DCwF has an operation
\begin{align*}
\mathsf{symm}&\colon\{\Gamma\colon\NCon\}\{A\colon\NTy\}\{t\colon\Tm(\Gamma,A^-)\}\{t'\colon\Tm(\Gamma,A)\}\\
    &\to \Tm(\Gamma,\Id(t,t')) \to \Tm(\Gamma,\Id(-t',-t))
\end{align*}

    \end{proposition}

    \begin{proof}

By the following construction in the DCwF syntax:

\includeCTFrame{directedTT/Id-symm}

    \end{proof}

This proof relies on the neutrality of $A$ in a very subtle, but critical way: in the definition of the type family $\mathsf{S}$, the variable term $\vv\colon\Tm(\Gamma\cExtend A,A[\pp_A])$ is negated, so that it is of type $A[\pp_A]^-$ and therefore able to stand as the first argument to $\Id$. But this is only possible if $\Gamma\cExtend A\colon\NCon$ because term-negation is only defined in neutral contexts. This reasoning will prove important for the style of reasoning we employ in \autoref{syntheticCT}, so we isolate it as a principle.

\begin{principle*}[Var Neg] \label{tag:RuleVarNegFormal}

For $\Gamma\colon\NCon$, the variable term $$\vv\colon\Tm(\Gamma\cExtend[+]A,A[\pp_{A}])$$  can only be negated (i.e. forming $-\vv$) if $A\colon\NTy\;\Gamma$ (and likewise for $\cExtend[-]$).
\end{principle*}
In \autoref{setUniverse}, we'll argue that there's \textit{no} way to construct this symmetry term (for arbitrary DCwFs\footnote{with some nontrivial amount of structure.}) if $A$ is not assumed to be neutral.

Before proceeding, it's worth explaining what is \enQuote{the DCwF syntax} mentioned in the proof above. This is where it becomes relevant that DCwFs are presented as generalized algebraic theories: as mentioned in the
introduction, \cite{kaposi2019constructing} proves that any GAT has an initial syntax model. Therefore, any construction done in the syntax model (such as the construction of $\mathsf{symm}$ above) can be interpreted into \textit{any} DCwF. This is why a syntactic construction was adequate to prove a claim about \textit{all} DCwFs in the foregoing proof. In the next section, our proofs will all be syntactic, and thereby apply to arbitrary DCwFs.

Let us make an important observation about the syntax of DCwFs. An important criterion for our theory is that hom-types can be \textit{iterated}, that is, our syntax allows for the formation of homs between homs, and homs between homs between homs, and so on. The iteration of identity types is, after all, how homotopy type theory is able to serve as a synthetic language for higher groupoids; and since hom types are iterable in the DCwF syntax, it is a synthetic language for higher categories. However, a given model may be \textit{truncated}, in that the higher structure may become trivial after a certain point. This is the case with the groupoid model: while its types do not all obey the \textit{uniqueness of identity proofs (UIP)} principle and are therefore not mere \textit{h-sets}, they do obey \enQuote{UIP, one level up}: in the groupoid model, identity proofs of identity proofs \textit{are} unique.

The same happens in the category model: in general, there may be terms $p\colon\Tm(\Gamma,\hom(t,t')^-)$ and $q\colon\Tm(\Gamma,\hom(t,t'))$ but no term of type $\hom(p,q)$. But if there is such a term, there is exactly one. So the (1-)category model, unsurprisingly, can only model 1-categories. But there's another sense in which the category model structure trivializes \enQuote{one level up}: all the hom-types are interpreted as discrete categories, which are necessarily \textit{groupoids}. So, while $\hom(t,t')$ is still a synthetic category, it's actually a synthetic \textit{setoid}. This is appropriate for doing synthetic 1-category theory: it makes sense that the hom-types are trivial \textit{as categories}: to do 1-category theoretic arguments, we wish to speak of \textit{identities} between parallel morphisms, not further category-theoretic structure.

We encapsulate DCwFs like this into a definition for further study.
\mkDefn{1-1-trunc-DCwF}
The numbering follows the well-known indexing of $(n,m)$-categories (see e.g. \cite[Defn. 8]{baez2007lecturesncategoriescohomology}) to refer to $\infty$-categories where all parallel $k$-morphisms are equal when $k>n$ and all $k$-morphisms are invertible for $k>m$. We could define $(n,m)$-truncated DCwFs for arbitrary $n$ and $m$ (for instance, the preorder model would be $(0,1)$-truncated, the groupoid model $(1,0)$-truncated, etc.), but that would take us too far afield. For the present work, we will work with $(1,1)$-truncated DCwFs, and develop the theory of synthetic 1-category theory (i.e. synthetic $(1,1)$-category theory) in that language. The practical consequence of working in the syntax of $(1,1)$-truncated DCwFs is that we only have one \enQuote{layer} of homs, and the type $\hom(t,t')$ itself is neutral, i.e. its homs are symmetric identity types.

To conclude this section, we state the $\Pi$-type former in PCwFs. This is just the appropriately-polarized analogue of \cite[Defn. 3.15]{hofmann1997syntax}, and is approximately the same rule for $\Pi$-types in \cite{licataHarper}.
\mkDefn{Pi}
The $\beta$-law is that $\code{app}\circ\code{lam}\equiv\iden$ and the $\eta$-law the other way around. The differences between polarized $\Pi$-types and the familiar $\Pi$-types of undirected type theory are pretty minimal when operating in neutral contexts: for instance, when instantiating to non-dependent functions $A \to B$, the application operator defined by $f\;\$\;t = (\app\;f)[\iden,_- t]$ has type
\[  \underlines\$\underlines \colon \Tm(\Gamma,A\to B) \to \Tm(\Gamma^-,A^-) \to \Tm(\Gamma,B). \]
If $\Gamma\colon\NCon$, then we can take terms from $\Tm(\Gamma^-,A)$, $\Tm(\Gamma,A[e])$ and $\Tm(\Gamma,A[e]^-)$ and use the negation operator and the $e^-$ substitution to get into $\Tm(\Gamma^-,A^-)$ for the purposes of applying functions. In the next section, we push this bureaucracy into the background and proceed informally.

\section{Synthetic Category Theory}\label{syntheticCT}
In this section, we work in an arbitrary $(1,1)$-truncated DCwF with polarized $\Pi$-types by only working in the syntax. In the main body of the text, we'll adopt an informal type theoretic style (inspired by \cite{hottbook}).
We assume that we're working in some neutral context $\Gamma$, though we don't explicitly reference $\Gamma$. We'll write $t\colon A$ to indicate $t\colon\Tm(\Gamma,A)$. In what follows, we'll use the letters $p,q,r,s,t,u,v,w,f,g$ to name terms (of various types) in $\Gamma$, whereas the letters $x,y,z$ will be the names of \textit{variables} obtained by extending $\Gamma$. We'll have to be careful to abide by the variable negation rule:
\begin{principle*}[Var Neg] \label{tag:RuleVarNegInformal}

An expression $e$ can only be negated if all the variables occurring in it are of neutral types.
\end{principle*}
We suppress the distinction between $\Gamma$ and $\Gamma^-$, since we can substitute back and forth with $e$ behind the scenes, as needed. Negative context extension will just behave like positive extension by a negative type: recall that $$\vv_{-,A}\colon\Tm((\Gamma\cExtend[-]A)^-,A[\pp_{-,A}^-]^-)$$ i.e. $$ \vv_{-,A}\colon\Tm(\Gamma^-\cExtend[+]A^-,A[\pp_{A^-}]^-) $$
so, if we're suppressing the distinction between $\Gamma$ and $\Gamma^-$, then this is just a variable $x$ of $A^-$. Accordingly, we'll apply functions like this:
\[  \infer{f(t)\colon B(t)}{f\colon \prod_{(x\colon A^-)} B(x) & t\colon A^-} \]
and form them like this:
\[  \infer{(\lambda (x\colon A^-) \to e)\colon \prod_{x\colon A^-}B(x).}{x\colon A^-\vdash e \colon B(x)} \]

Finally, here's our principle of directed path induction:
\begin{principle*}[Directed Path Induction] \label{tag:RuleDirJ}

For every $t\colon A^-$, if $M(x,y)$ is a type family depending on $x\colon A$ and $y\colon\hom(t,x)$, then, for each $$m\colon M(-t,\refl_t),$$we get an $$\ind_M(m,x,y)\colon M(x,y)$$ for all $x,y$.
\end{principle*}
So, for instance, the construction of symmetry above (the proof of \propRef{Id-symm}) would be expressed informally as follows: given a neutral type $A$ and a term $t\colon A^-$, define a type family over $x\colon A,y\colon Id(t,x)$ by
\[ S(x,y) = \Id(-x,-t) \]
We have not violated \ruleRef{RuleVarNegInformal} because $x\colon A$ and $A$ is neutral. We have a term of type $S(-t,\refl_t)$, i.e. $\Id(t,-t)$, namely $\refl_t$. So therefore we get $S(x,y)$ for arbitrary $x,y$. If we have a particular $t'\colon A$ and $p\colon\Id(t,t')$, we can put $$ \mathsf{symm}\;p = \ind_S(\refl_t,t',p).$$ Again, we emphasize that it is \ruleRef{RuleVarNegInformal} which prevents this argument from working for non-neutral types, as desired. Below, we are more casual with our application of directed path induction (e.g. not defining the type family explicitly) in cases where \ruleRef{RuleVarNegInformal} is not a concern.

With that, we can proceed to the informal constructions. Along the way, the explicit constructions in the DCwF syntax are carried out in the accompanying figures.

\subsection{Composition of Homs}
\textit{Formal development: \autoref{tag:HomComposition} }

As mentioned, a type $A$ in directed type theory is supposed to be a \textit{synthetic category}. The terms $t:A$ represent objects, and the terms $p\colon\hom(t,t')$ represent morphisms. For this to truly be category theory, however, we must be able to compose morphisms. We'll write composition in diagrammatic order: given $t,u\colon A^-$ and $v'\colon A$, we should be able to compose $p\colon\hom(t,-u)$ with $q\colon\hom(u,v')$ to get $p\cdot q\colon\hom(t,v')$. We do this by directed path induction on $q$, by putting
\[    p\cdot\refl_{u} = p. \]
The $\refl$ terms serve as the identity morphisms of the category: by the above, we know that $p\cdot\refl_{u} \equiv p$, and thus $\refl_p\colon \Id(p\cdot\refl_{u},p)$. As for the other unit law, we must again use directed path induction: since
  \[  \refl_{u}\cdot \refl_{u} \equiv \refl_{u}, \]
we have that $\refl_{\refl_{u}}\colon\Id(\refl_{u}\cdot\refl_{u},\refl_{u})$, and, by induction we get a term
\[ \mathsf{r{-}unit}\;q = \ind(\refl_{\refl_u},v,q)\colon \Id(\refl_{u}\cdot q,q)  \]
for each $q\colon\hom(u,v')$.

Finally, we get that the composition operation is associative. Given $t,u,v\colon A^-$ and $w'\colon A$ as well as $p\colon\hom(t,-u)$, $q\colon\hom(u,-v)$, and $r\colon\hom(v,w')$, we construct
\[ \mathsf{assoc}\;p\;q\;r\colon\Id(p\cdot(q\cdot r), (p\cdot q)\cdot r) \]
by directed path induction on $r$. If $r=\refl_{v}$, then $q\cdot r\equiv q$ and $(p\cdot q)\cdot r\equiv p\cdot q$. Thus, we have
\[ \refl_{p\cdot q} \colon \Id(p\cdot(q\cdot \refl_{v}), (p\cdot q)\cdot \refl_{v}) \]
and then the induction carries through, and we get $\mathsf{assoc}\;p\;q\;r$ as desired.

\begin{FigCT}{categoryModel/HomComp}{HomComposition}

\caption{Composition of Homs}
\end{FigCT}

\subsection{Synthetic Functors}
\textit{Formal development: \autoref{tag:map}, \autoref{tag:mapCalc}, \autoref{tag:FuncComposition}, \autoref{tag:FuncCompCalc} }

If types $A,B$ are synthetic categories, it should come as no surprise that terms $f\colon A\to B$ are synthetic \textit{functors}. The object part is given by the usual function application, but the variances are somewhat mixed: if $t\colon A^-$, then we can say $f(t)\colon B$. However, we can still apply $f$ to a term $t'\colon A$, we just have to put a minus on $t'$, i.e. $f(-t')$.

Unlike usual (\enQuote{analytic}) category theory, we don't have to explicitly define the morphism part of a functor; any term of type $A\to B$ we can write down will come with a morphism part for free. To obtain this morphism part, again we use directed path induction: given an $f\colon A\to B$ and some $t\colon A^-$, we can define a $B$-morphism
\[\map\;f\;p \colon \hom(-f(t),f(-t')) \]
for every $t'\colon A^-$ and $p\colon\hom(t,t')$ by putting
\[
    \map\;f\;\refl_t = \refl_{-f(t)} \quad\colon\hom(-f(t),f(t)).
\]
By definition, this operation preserves identities (sending $\refl$ to $\refl$), and respects composition: if we have $t,u\colon A^-$ and $p\colon\hom(t,-u)$, then, since $\map\;f\;\refl_{u}\equiv \refl_{-f(u)}$ and $p\cdot\refl_{u}\equiv p$ and $(\map\;f\;p)\cdot\refl_{-f(u)}\equiv \map\;f\;p$, we have
\[
\refl_{(\map\;f\;p)}\colon\Id(\map\;f\;(p\cdot\refl_{u}), (\map\;f\;p)\cdot(\map\;f\;\refl_t)).
\]
By induction, we get an identity between $\map\;f\;(p\cdot q)$ and $(\map\;f\;p)\cdot(\map\;f\;q)$ for arbitrary $q$.

Let us also note that functors are also composable: given $f\colon A\to B$ and $g\colon B\to C$, we get the usual
\[ g\circ f = \lambda (x : A^-) \to g(-f(x)). \]
Of course, we can prove $\map\;(g\circ f)\;p$ equal to $\map\;g\;(\map\;f\;p)$ by directed path induction, using the following observations:
\begin{align*}
   \map\;f\;\refl_t &\equiv \refl_{-f(t)}\\
   \map\;g\;\refl_{-f(t)} &\equiv \refl_{-g(-f(t))}\\
   \map\;(g\circ f)\;\refl_t &\equiv \refl_{-(g\circ f)(t)}.
\end{align*}

\begin{FigCT}{categoryModel/map}{map}

\caption{Morphism part of Functors}
\end{FigCT}

\begin{figure}

\begin{align*}
    &\hom(-(f\;\$\;t), -(-(f\;\$\;t)))\\
    &\equiv \hom(-(f\;\$\;t), f\;\$\;t )\\
    &= \hom(-(f\;\$\;t), (\app\;f)[\iden ,_- t]) \\
    &\equiv \hom(-(f\;\$\;t), (\app\;f)[ (e\cExtend A) \circ \pair{\iden ,_+ -t[e]}]) \tag{*}\\
    &\equiv \hom(-(f\;\$\;t), (\app\;f)[(e\cExtend A)][\iden ,_+ -t[e]]) \\
    &\equiv \hom(-(f\;\$\;t)[\pp_{A[e]}, (\app\;f)[(e\cExtend A)])[\iden ,_+ -t[e]] \\
    &= \mathsf{MAP}[\iden ,_+ -t[e]]
\end{align*}
(*) $\pair{\iden ,_- t} \equiv (e\cExtend A)\circ\pair{\iden ,_+ -t[e]}$ by \defnRef{NPCwF} since $\Gamma$ neutral.
\caption{Calculation that $\refl_{-(f\;\$\;t)}\colon \mathsf{MAP}[-t[e]]$ in \autoref{tag:map} } \label{tag:mapCalc}
\end{figure}

\begin{FigCT}{categoryModel/FuncComp}{FuncComposition}

\caption{Composition of Functions}
\end{FigCT}

\begin{figure}
\raggedright
The diagram
\[
\begin{tikzcd}[column sep=large]
    \Gamma\cExtend[-]A
        \arrow[from=1-1, to=1-2, "\pair{\pp_{-,A} ,_+ \app\;f}" {yshift=5pt}] \arrow[bend right=10, shift right=4, swap]{rrr}{\pp_{-,A}}  & \Gamma\cExtend[+]B \arrow[swap]{r}{e\cExtend B[e^-]} \arrow[bend left=8, shift left=4]{rr}{\pp_B} & \Gamma\cExtend[-] B[e^-] \arrow{r}{\pp_{-,B[e^-]}} & \Gamma
\end{tikzcd}
\]
commutes, since $p_{-,B[e^-]}\circ (e\cExtend B[e^-])\equiv\pp_B$ by \defnRef{NPCwF} and $\pp_B\circ\pair{\pp_{-,A},_+\app\;f} \equiv \pp_{-,A}$ by \autoref{tag:CwFDefn}. Thus, since $\app\;g\colon\Tm(\Gamma\cExtend[-]B[e^-],C[\pp_{-,B[e^-]}])$,
\[
(\app\;g)[e\cExtend B[e^-]][\pp_{-,A},_+\app\;f] \colon \Tm(\Gamma\cExtend[-]A,C[\pp_{-,A}]).
\]
\caption{Auxiliary calculations for \autoref{tag:FuncComposition} } \label{tag:FuncCompCalc}
\end{figure}

\section{Further observations about the Category Model}\label{setUniverse}

As the previous section showed, the syntax of 1-truncated Directed CwFs provides a nice setting for some very basic constructions in synthetic category theory. However, further expansion of the DCwF syntax is needed to be able to capture the full range of constructions in category theory. In this section, we'll observe some constructions that can be made (and some equivalences that hold) in the category model, which require further study to be internalized into the DCwF syntax.

Probably the most significant omission from the synthetic category theory of the previous section is natural transformations. There are some natural transformations expressible in the theory as written, because natural transformations are, as we might hope and expect, homs between functors. That is, the type $A\to B$ is a synthetic category: given $f\colon\Tm(\Gamma,(A\to B)^-)$ and $g\colon\Tm(\Gamma,A\to B)$, we can form the type $\hom(f,g)$. In the category model, these are interpreted as natural transformations from $f$ to $g$, but all dependent over the context $\Gamma$. Both $f$ and $g$ sends objects $\gamma\colon\abs{\Gamma}$ to functors from $A(\gamma)$ to $B(\gamma)$; if $\alpha\colon\Tm(\Gamma,\hom(f,g))$, then $\alpha$ sends $\gamma$ to a natural transformation $f(\gamma)\to g(\gamma)$. On morphisms $\gamma_{01}\colon\Gamma\:[\gamma_0,\gamma_1]$, $\alpha(\gamma_{01})$ is a witness to the fact that
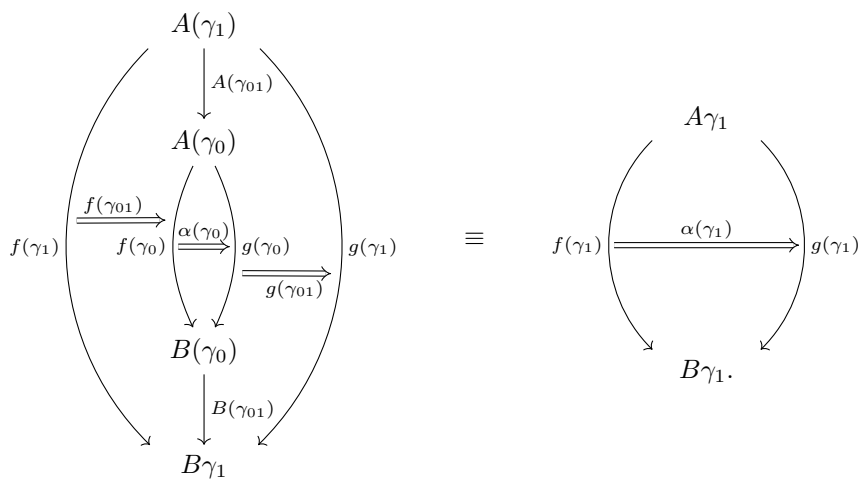
\[
    \begin{tikzcd}[sep=large]
        A(\gamma_1)
            \arrow[swap,bend right=45,shift right=5,""{name=LL}]{dddd}{f(\gamma_1)}
            \arrow[bend left=45,shift left=5,""{name=RR}]{dddd}{g(\gamma_1)}
            \arrow{d}{A(\gamma_{01})}\\
        A(\gamma_0)
            \arrow[swap, bend right=25,""{name=Ll}]{dd}{f(\gamma_0)}
            \arrow[bend left=25,""{name=Rr}]{dd}{g(\gamma_0)} \\
            \ar[from=LL, to=Ll, shift left=10pt,"f(\gamma_{01})", Rightarrow, shorten <= 4pt, shorten >= -2pt]
            \ar[from=Rr, to=RR, swap, shift right=10pt, "g(\gamma_{01})", Rightarrow, shorten <= -2pt, shorten >= 4pt]
            \ar[from=Ll, to=Rr, "\alpha(\gamma_{0})", Rightarrow, shorten <= 2pt, shorten >= 2pt]
        \\
        B(\gamma_0)\arrow{d}{B(\gamma_{01})} \\
        B \gamma_1
    \end{tikzcd}
    \qquad
    \equiv
    \qquad
    \begin{tikzcd}
        A \gamma_1
            \arrow[swap,bend right=45,shift right=5,""{name=LL}]{dddd}{f(\gamma_1)}
            \arrow[bend left=45,shift left=5,""{name=RR}]{dddd}{g(\gamma_1)}\\
            \ar[from=LL, to=RR, "\alpha(\gamma_{1})", Rightarrow, shorten <= 2pt, shorten >= 2pt]
        \\
        \\
        \\
        B \gamma_1.
    \end{tikzcd}
\]
Presently, the only such natural transformations expressible in the DCwF syntax are identities (e.g $\refl_f\colon\Tm(\Gamma,\hom(f,-f))$). There is not a way to work with the actual components of a natural transformation, or to define a natural transformation \textit{by} its components. What we would like to be able to do is write terms of type
\[ \prod_{x\colon A^-} \hom(-f(x),g(x)) \]
and then prove that they are automatically natural by directed path induction, since $\alpha_t\cdot \map\;g\;\refl_t \equiv \map\;f\;\refl_t \cdot \alpha_t$. The issue is that this type violates \ruleRef{RuleVarNegInformal}: we're not allowed to write $-f(x)$ for a variable $x$ without knowing that $A$ is neutral. If $A$ \textit{is} neutral, then we're only capturing natural transformations between functors whose domain is a groupoid, which is a significant restriction. Thus, we have two important open questions: how to make natural transformations definable component-wise in the syntax (ideally using $\Pi$-types), and how to express in the syntax that the type of such transformations is equivalent to the hom-type between the two functors.

Another important feature we add are \textit{universes}. The category model comes equipped with several type universes, most significantly the universe of sets. More precisely, we can regard the category $\SET$ as a closed type in the category model. The operation $\El\colon\Tm(\Empty,\SET) \to \Ty\;\Empty$ takes a set $X$ and views it as a discrete category. We can then define
\[  \code{Hom-to-func} \colon \{X\colon\Tm(\Empty,\SET^-)\}\{Y\colon\Tm(\Empty,\SET)\} \to \Tm(\Empty,\hom(X,Y)) \to \Tm(\Empty,\El(X)\to\El(Y)) \]
by directed path induction: $\code{Hom-to-func}\;\refl_X$ should be the identity function $\mathsf{lam}\;\vv\colon\Tm(\Empty,\El(X)\to\El(X))$. We can use this to state the following principle.
\begin{principle*}[External Directed Univalence] \label{tag:dirUnivalence}

\code{Hom-to-func} is a bijection.
\end{principle*}
Really, \code{Hom-to-func} is the identity function, since terms of hom-types in the empty context are just the morphisms of the category, the morphisms of $\SET$ are functions, and a functor between discrete categories is just a function between their objects. Spelling out the category model semantics, we see that every function is sent to itself. Sufficiently internalized, this principle of Directed Univalence serves as the directed analogue of Hofmann and Streicher's \textit{universe extensionality} axiom \cite[Section 5.4]{groupoidModel}. Further work is required to better develop the theory of isomorphisms in the synthetic category theory, and to compare this principle of directed univalence to existing ones (e.g. \cite{kudasov2023formalizing}).

Let us conclude by observing that the existence of a universe allows for metatheoretic reasoning as well, specifically negative proofs about what \textit{cannot} be done in the syntax. We can view the directed universe $\SET$ as a source of nontrivial directedness: if we affirm \autoref{tag:dirUnivalence}, then $\SET$ cannot possibly be a neutral type. We show that, in DCwFs equipped with a directed univalent universe $\SET$, hom-types must be asymmetric in general. That is, we cannot construct a term $\mathsf{symm}$ like in \propRef{Id-symm} for non-neutral types in the syntax DCwF+$\SET$ (the initial model of the GAT of DCwFs with a universe $\SET$). We do so the same way Hofmann and Streicher \cite{groupoidModel} proved that ordinary Martin-L{\"o}f Type Theory couldn't prove the Uniqueness of Identity Proofs: by countermodel. Hofmann and Streicher's countermodel was the groupoid model, and, of course, ours is the category model.
    \begin{proposition}\label{state:cant-prove-symm}
There cannot be an operation
\begin{align*}
\mathsf{symm'}&\colon\{\Gamma\colon\NCon\}\{A\colon\Ty\}\{t\colon\Tm(\Gamma,A^-)\}\{t'\colon\Tm(\Gamma,A)\}\\
    &\to \Tm(\Gamma,\hom(t,t')) \to \Tm(\Gamma,\hom(-t',-t))
\end{align*}
definable in the syntax of DCwFs+$\SET$.

    \end{proposition}

    \begin{proof}
If the syntax model of DCwF+$\SET$ had such an operation $\mathsf{symm'}$, then, by initiality, so too would every DCwF with $\SET$, in particular the category model. But then for any $X\colon\Tm(\Empty,\SET^-)$ and $Y\colon\Tm(\Empty,\SET)$ and $f\colon\Tm(\Empty,\hom(X,Y))$, we would obtain $\mathsf{symm'}\;f\colon\Tm(\Empty,\hom(Y,X))$. But this is absurd, because the function $?\colon\emptyset\to\mathbbm{1}$ is a term of type $\El(\emptyset)\to\El(\mathbbm 1)$ in the category model, and, by \autoref{tag:dirUnivalence}, corresponds to a term of type $\hom(\emptyset,\mathbbm{1})$, but there cannot be any terms of $\hom(\mathbbm 1,\emptyset)$, because the set of terms of this type is in bijection with the set of functors  $\El(\mathbbm 1)$ to $\El(\emptyset)$, of which there are none.

    \end{proof}

Basically the same argument will show the \textit{uniqueness of homs} principle---that for any hom terms $p\colon\Tm(\Gamma,\hom(t,t')^-)$ and $q\colon\Tm(\Gamma,\hom(t,t'))$, there is a witness of $\Id(p,q)$---is violated in the category model (a counterexample being $\SET$-homs from the two-element set to itself), and therefore not provable in the syntax of DCwFs+$\SET$. So we can conclude that the difference between $(1,1)$-truncated DCwFs, $(1,0)$-truncated DCwFs, and $(0,1)$-truncated CwFs is reflected internally in the syntax.

\section{Conclusion and Future Work}\label{conclusion}
We have laid the laid the foundation for a generalized algebraic theory of directed types, and began to conduct synthetic category theory in that setting. Our semantics-forward approach was to study the category model first, and extract its key features into a series of abstract definitions---the GATs of polarized CwFs, neutral-polar CwFs, directed CwFs, $(1,1)$-truncated directed CwFs, and directed CwFs with features like polarized dependent types and a directed univalent set universe. Working within the directed type theory of these models, we found that it was possible to work informally with the powerful directed path induction principle to make basic constructions in category theory, with our careful discipline about variable negation preventing the directed type theory from collapsing into undirected type theory.

Much remains to be done. The category theory of \autoref{syntheticCT} serves as a proof-of-concept, but needs to be fleshed out into a full theory. As mentioned, work is needed to articulate natural transformations in the theory; our current investigations concern possible generalization of $\Pi$-types to di-variant \textit{end types} which address some of the above-mentioned variance issues with natural transformations. For reasons of space, we omitted dependent sum types from the theory. But with them added, much of basic category theory should be expressible in this language, such as isomorphisms, (co)slice categories, (co)limits, exponentials, and perhaps some basic topos theory. Better development of the category of sets should put representability and some properties of presheaf categories into reach, though internal statement and proof of the Yoneda Lemma will likely rely on the resolution of the above-mentioned dilemma regarding natural transformations. We also leave it to future work to study whether this theory can capture higher category theory by weakening or dropping the assumption of 1-truncation, and, if so, how it compares to existing synthetic higher category frameworks, such as \cite{riehl2017}.

There are further avenues for developing the type theory of DCwFs. Two important metatheoretic results about the syntax of DCwFs currently being pursued are \textit{canonicity} and \textit{normalization}. Moreover, we would like to verify the correctness of these results by formalizing them in a computer proof assistant. A further goal would be to implement the syntax of DCwFs as a computer proof language itself, hopefully with syntax nearly as convenient as the constructions of \autoref{syntheticCT}, and formalize larger swaths of category theory in it.

As mentioned in the introduction, a motivation for the present work's focus on generalized algebraic theories is the possibility of expressing it in a \textit{second-order} generalized algebraic theory, following \cite{bocquet2022external, Uemura2023, uemura2021abstract, altenkirch2024internal, kaposi2024second}. Since our notions of PCwFs and NPCwFs include explicit operations on contexts (as seen in the substructural character of the \ruleRef{RuleVarNegInformal} rule), it's clear that either an extension to the SOGAT signature language of \cite{kaposi2024second}, and/or a partial internalization of the first-order theory into the second-order theory---{\`a} la \cite{altenkirch2024internal}---will be necessary. A related question is whether the $(\underlines)^-$ operation (or a neutralization operation interpreted in the category model by core groupoids) on types and contexts can be viewed as a modality in the sense of \cite{gratzer2020multimodal}. Also, as mentioned, this work is oriented towards \textit{higher observational type theory}, and further study of the observational equivalences of this theory (e.g. a characterization of the $\hom$-types of $\Pi$-types) is needed.

Finally, the present framework provides a setting for studying \textit{directed higher-inductive types}---inductively-defined types with both term constructors and \textit{hom constructors}. Some simple examples would be the directed interval (as is studied in \cite{weaverLicata, riehl2017}) and a directed analogue of the \textit{circle} type \cite[Section 6.2]{hottbook}. These examples are modelled by the category model, and can therefore be soundly added to the present theory. Higher examples (such as directed versions of higher tori and spheres) would require more careful metatheoretic work to justify, but could perhaps lead to a number of interesting considerations.

\bibliographystyle{plainurl}
\bibliography{biblio}

\appendix

\section{Additional Calculations}\label{appendixCalcs}

The definition of the $\Pi$-type former given in \autoref{tag:CatModelPi} is repeated in \autoref{tag:CatModelPiFull}, plus the definitions for lambda abstraction and application.

These definitions rely on the following calculations.
\begin{itemize}
\item A naturality calculation for the morphism part of \code{lam t'}:
\begin{equation}\label{tag:PiTypeNaturality}
\begin{aligned}
    & t'(\gamma_1,x_1)\circ B(\iden_{\gamma_1},x_1)(t'(\gamma_{01},\iden_{A\;\gamma_{01}\;a_1}))\\
    &\equiv t'(\iden_{\gamma_1}\circ \gamma_{01}, (A\;\gamma_{01}\;x_1)\circ \iden_{A\;\gamma_{01}\;a_1})\\
    &\equiv t'(\gamma_{01}, A\;\gamma_{01}\;x_1)\\
    &\equiv t'(\gamma_{01}\circ\iden_{\gamma_0}, (A\;\iden_{\gamma_0}\;\iden_{A\;\gamma_{01}\;a_1'}) \circ A\;\gamma_{01}\;x_1)\\
    &\equiv t'(\gamma_{01},\iden_{A\;\gamma_{01}\;a_1'}) \circ B(\gamma_{01},\iden_{A\;\gamma_{01}\;a_1'})(t'(\iden_{\gamma_0},A\;\gamma_{01}\;x_1))
\end{aligned}
\end{equation}
The first and the last equations are the functoriality conditions of $t'$. The middle equations are the category laws for $\Gamma$, $A(\gamma_0)$, and $A(\gamma_1)$, as well as the functoriality of $A$.

\item To see that the morphism part of \code{app f} is well-typed, observe that
\begin{equation}\label{tag:appMorphismCalc1}
\begin{aligned}
&(f\;\gamma_{01}) \colon (\Pi(A,B)\;\gamma_1)[\;\Pi(A,B)\;\gamma_{01}\;(f\;\gamma_0),\;f(\gamma_1)\;], \\
&\text{i.e.}\\
&(f\;\gamma_{01}) \colon \mathsf{Transform}\;\{\gamma=\gamma_1\}\;(B(\gamma_{01},\iden)\;(f\;\gamma_0), f(\gamma_1))\\
&\text{and therefore,}\\
&\mathsf{component}(f\;\gamma_{01})\;a_1\\ &\colon B(\gamma_1,a_1)[ B(\gamma_{01},\iden)\;(f\;\gamma_0\;(A\;\gamma_{01}\;a_1)), f\;\gamma_1\;a_1];
\end{aligned}
\end{equation}
and also that
\begin{equation}\label{tag:appMorphismCalc2}
\begin{aligned}
&a_{01}\colon (A\;\gamma_0)[ a_0, A\;\gamma_{01}\;a_1 ], \\
&\text{and thus, since }(f\;\gamma_0)\colon\Tm(A(\gamma_0),B_{\gamma_0}),\\
&f\;\gamma_0\;a_{01}\\ &\colon B(\gamma_0,A\;\gamma_{01}\;a_1)[ B(\iden_{\gamma_0},a_{01})\;(f\;\gamma_0\;a_0), f\;\gamma_0\;(A\;\gamma_{01}\;a_1)]\\
&\text{and therefore}\\
&B(\gamma_{01},\iden)\;(f\;\gamma_0\;a_{01})\\
&\colon B(\gamma_1,a_1)[ B(\gamma_{01},a_{01})\;(f\;\gamma_0\;a_0), B(\gamma_{01},\iden)\;(f\;\gamma_0\;(A\;\gamma_{01}\;a_1)) ].
\end{aligned}
\end{equation}
and thus we can conclude that $\mathsf{component}(f\;\gamma_{01})\;a_1$ can be composed with $B(\gamma_{01},\iden)\;(f\;\gamma_0\;a_{01})$, giving a term of the appropriate type.

\end{itemize}

\begin{FigCT*}[p]{categoryModel/Pi-lamapp}{CatModelPiFull}

\caption{Complete semantics of $\Pi$-types in the category model}
\end{FigCT*}

\end{document}